\documentclass[12pt]{article}

\usepackage{epsfig,amsmath,amssymb,graphics,verbatim,amsfonts,subfigure,psfrag}
\usepackage[letterpaper,margin=1in]{geometry}

\let\oldbibliography\thebibliography
\renewcommand{\thebibliography}[1]{%
  \oldbibliography{#1}%
  \setlength{\itemsep}{0.5mm}%
}

\def\R{\mathbb{R}}

\def\I{\infty}

\def\txtd{{\textnormal{d}}}

\def\txtD{{\textnormal{D}}}

\newcommand{\be}{\begin{equation}}
\newcommand{\ee}{\end{equation}}
\newcommand{\bea}{\begin{eqnarray}}
\newcommand{\eea}{\end{eqnarray}}
\newcommand{\beann}{\begin{eqnarray*}}
\newcommand{\eeann}{\end{eqnarray*}}
\newcommand{\benn}{\begin{equation*}}
\newcommand{\eenn}{\end{equation*}}

\def\ra{\rightarrow}
\def\I{\infty}

\newcommand{\cH}{{\mathcal H}}  
\newcommand{\cL}{{\mathcal L}}  
\newcommand{\cM}{{\mathcal M}}  
\newcommand{\cN}{{\mathcal N}}  
\newcommand{\cT}{{\mathcal T}}  

\begin{document}

\author{Christian Kuehn\thanks{Institute for Analysis and Scientific Computing, 
Vienna University of Technology, 1040 Vienna, Austria}}
\title{Efficient Gluing of Numerical Continuation and a Multiple Solution Method for Elliptic PDEs}
\maketitle
\begin{abstract}
Numerical continuation calculations for ordinary differential equations (ODEs) are, by now, an 
established tool for bifurcation analysis in dynamical systems theory as well as across almost
all natural and engineering sciences. Although several excellent standard software packages are available for ODEs,
there are - for good reasons - no standard numerical continuation toolboxes available for
partial differential equations (PDEs), which cover a broad range of different classes of PDEs automatically. 
A natural approach to this problem is to look for efficient gluing 
computation approaches, with independent components developed by researchers in numerical 
analysis, dynamical systems, scientific computing and mathematical modelling. In this paper, we shall
study several elliptic PDEs (Lane-Emden-Fowler, Lane-Emden-Fowler with microscopic force, 
Caginalp) via the numerical continuation software pde2path and develop a gluing component to
determine a set of starting solutions for the continuation by exploiting the variational 
structures of the PDEs. In particular, we solve the initialization problem of numerical continuation for PDEs
via a minimax algorithm to find multiple unstable solution. Furthermore, for the 
Caginalp system, we illustrate the efficient gluing link of pde2path to the underlying mesh
generation and the FEM MatLab pdetoolbox. Even though the approach works efficiently due to the high-level
programming language and without developing
any new algorithms, we still obtain interesting bifurcation diagrams and directly
applicable conclusions about the three elliptic PDEs we study, in particular with respect to symmetry-breaking. 
In particular, we show for a modified Lane-Emden-Fowler equation with an asymmetric microscopic
force, how a fully connected bifurcation diagram splits up into C-shaped isolas on which localized
pattern deformation appears towards two different regimes.
We conclude with a section on future software development issues that would be helpful to be addressed to 
simplify interfaces 
to allow for more efficient, time-saving, gluing computation for dynamical systems analysis of PDEs 
in the near future.
\end{abstract}


\section{Introduction}
\label{sec:intro}

One of the most common problems encountered in the analysis of differential equations arising
in applications is to give a qualitative description of the influence 
of parameters $\mu\in\R^p$ on the system dynamics. For ordinary differential equations (ODEs) with
phase space variables $x\in \R^N$, one approach to this problem is to use a standard 
numerical integration method in combination with Monte-Carlo 
simulation. For obtaining a basic idea about the system dynamics and exploration of different regimes,
this approach is an indispensable tool in modern scientific computing. However, to completely map out
the different dynamical regimes in parameter space, this approach requires 
a systematic simulation over grids of initial 
values as well as over grids in parameter space. Even though there are techniques available to 
alleviate some of the computational effort when $p$ and $N$ are large (formally we write this as 
$p,N\gg1$), such as multi-level Monte Carlo (MC) \cite{Heinrich} and quasi-MC \cite{Caflisch}, 
the computational effort can still be substantial. In addition, manual inspection of trajectories 
is usually required to detect important invariant sets, metastable states, bifurcations or locally 
invariant manifolds, which characterize the dynamics. Therefore, direct simulation can be a
computation- and labor-intensive approach to determine dynamical structures in nonlinear systems\footnote{
The MC approach with manual inspection of orbits is actually used very successfully - 
particularly in many application areas - via a method that could be called gsMC (graduate student 
Monte-Carlo).}. Indeed, this problem is well-known for some time as discussed in 
\cite[p.1-2]{TuckermanBarkley}:
\textit{``However, these [time-stepping] techniques rely on waiting out slow exponential decay,
or on arduous binary searches, and make highly inefficient use of both machine and human resources.''}
\medskip

As an alternative, numerical continuation methods are frequently used across all the natural 
sciences and engineering to systematically explore parameter space for important dynamical 
structures \cite{KrauskopfOsingaGalan-Vioque,KuepperSeydelTroger,RooseDeDierSpence,Seydel}. 
This approach has been pioneered in the late 1970s, particularly by Keller 
\cite{Keller,Keller3,Keller2}. The main idea is to use a predictor-corrector scheme in
combination with singularity detection to track parametrized curves (or higher-dimensional 
manifolds) of solutions to nonlinear algebraic systems arising, {e.g.}, from the steady
state formulation of the differential equation. The theory of numerical continuation is 
well-established by now \cite{AllgowerGeorg,DoedelKellerKernevezI,DoedelKellerKernevezII,Govaerts}. For 
ODEs, it has been used to compute, {e.g.}, bifurcation diagrams \cite{Doedel1,Kuznetsov}, stable/unstable 
manifolds \cite{KrauskopfOsinga1,KrauskopfSurvey}, global homoclinic/heteroclinic orbits 
\cite{Beyn2,ChampneysKuznetsovSandstede}, canards and canard values in fast-slow systems 
\cite{DesrochesKrauskopfOsinga2,KuehnCanLya}, invariant tori \cite{SchilderOsingaVogt}, 
spectra \cite{RademacherSandstedeScheel} and higher-dimensional solution manifolds \cite{Henderson1}. 
There are also several highly successful software packages available for ODE continuation problems 
such as \texttt{AUTO} \cite{Doedel_AUTO2007}, \texttt{XPP} \cite{XPP}, \texttt{MatCont} \cite{DhoogeGovaertsKuznetsov} or \texttt{PyDSTool} \cite{ClewleySherwoodLaMarGuckenheimer}. A key reason, why one may provide excellent
self-contained tools for ODEs, as evidenced by the structure of \texttt{AUTO} \cite{Doedel_AUTO2007}, is
that many questions may be formulated in terms of classical two-point boundary value problems (BVPs)
\cite{AscherMattheijRussell}. Although providing a suitable problem reformulation is not a trivial
task, there is substantial practical evidence that the strategy to continue two-point boundary-value 
problems is very successful\footnote{One may borrow terminology from physics and describe the 
situation regarding two-point BVPs as the main ``universality class'' in dynamical systems
computation for ODEs.}.\medskip

Furthermore, similar approaches as for ODEs have also been implemented for delay differential
equations (DDEs) \cite{EngelborghsLuzyaninaSamaey,SzalaiStepanHogan} and have shown promise to also
work for stochastic differential equations (SDEs) \cite{BarkleyKevrekidisStuart,KuehnSDEcont1}.
Hence, it is quite natural to also consider partial differential equations and numerical continuation.
This approach has been carried out, {e.g.}, for the Ginzburg-Landau equation 
\cite{SchloemerAvitabileVanroose}, excitable media \cite{Barkley1,KneesTuckermanBarkley}, 
thin film equations \cite{BeltrameThiele}, the Swift-Hohenberg
equation \cite{Avitabileetal,LloydSandstedeAvitabileChampneys}, two-dimensional ocean models 
\cite{KatsmanDijkstraSchmeits}, the Schnakenberg system \cite{UeckerWetzel}, among certainly many
other systems. However, there is no uniform approach via a single software package and/or 
continuation toolbox available.\medskip

In fact, it is unlikely that we are ever going to reach a similar status in terms
of stand-alone software, as {e.g.}~\texttt{AUTO} and \texttt{MatCont} for ODEs, in the case 
of numerical continuation for PDEs, not
even for the stationary case. The main reason is that for each class of PDEs, each type of spatial 
domain and each type of boundary conditions, there is an vast number of different spatial discretization 
methods available such as finite element methods (FEM) \cite{BrennerScott}, finite-volume methods (FVM) 
\cite{LeVeque1}, finite-difference methods (FDM) \cite{Thomas}, spectral methods \cite{GottliebOrszag}, 
and many others \cite{ArnoldBrezziCockburnMarini,BanerjeeButterfield,CottrellHughesBazilevs,Peskin}. 
Furthermore, problems arise on the linear algebra side, when trying to solve large linear systems arising 
from various discretizations \cite{Bindeletal,Dijkstraetal}. Given a PDE from applications, the natural question
is, what method to use within a numerical continuation framework? The answer clearly is: 
It depends! Indeed, the various classes of PDEs are so different from a numerical analysis viewpoint,
that all we may hope for is to have several software libraries available, one for each class and 
each method. If we want to combine numerical continuation with numerical PDE techniques, there are 
two natural approaches. One may just implement the discretization scheme for the PDE ``by-hand'' for
each problem and then add on top of this scheme ``by-hand'' a numerical continuation approach. This
has been used successfully in the PDE continuation examples mentioned above. However, the potential 
drawback is cost: namely computational cost due to potentially sub-optimal implementation of each 
component and labor cost due to the programming from first principles. The last issue was nicely
formulated - in a slightly different context - by Bangerth and Heister \cite{BangerthHeister}; just 
one important quote from their opinion 
piece summarizes the situation well: \textit{``If every graduate student writes the code for a new 
discretization from scratch, we will be stuck forever solving `toy problems' (like the proverbial
`Laplace equation on the square'). Unfortunately, this is no longer sufficient to convince our 
colleagues in the applied sciences [...] that the new method is also applicable to their vastly
more complex problems.''}\footnote{One may view the described situation as another incarnation of 
gsMC.}\medskip 

It is one main aspect - among several - of this paper to show that the combination of numerical continuation,
PDE discretization, scientific software packages and suitable gluing codes, is a
convincing strategy to solve relatively complex problems at low computational and low labor costs. 
Of course, one may argue that \textit{in mathematical theory} there is no problem with this well-known 
idea, it does work. However, \textit{in numerical practice}, the situation is quite different. Even
though software for gluing computations within a certain framework, such as \texttt{Trilinos} \cite{Herouxetal}, 
are available, there are still major practical problems to be solved for dynamical systems analysis of 
PDEs. However, to illustrate that we are right now in the situation, where aiming to glue several algorithms and
software packages in an efficient way, can be successful, we study here various elliptic PDEs:
the Lane-Emden-Fowler equation \cite{Chandrasekhar,Wong} with a linear term, the Lane-Emden-Fowler 
equation with a localized force and the Caginalp 
system \cite{Caginalp,CaginalpFife}. The main computational framework is set up in \texttt{MatLab} 
as a high-level gluing language with a
focus on the recently developed numerical continuation software \texttt{pde2path} \cite{UeckerWetzelRademacher,DohnalRademacherUeckerWetzel}.\medskip 

It is well-known by users of numerical continuation software that, 
during the practical implementation process, finding starting solutions is a key obstacle as 
pointed out nicely by Sandstede and Lloyd
\cite{SandstedeLloyd}: \textit{``Strategies 
for finding good initial guesses depend strongly on the underlying specific problem: often, 
finding good starting data is the main obstacle for using [numerical continuation] successfully''}\footnote{
The original quote contains [AUTO] instead of [numerical continuation] but it is clear
that the problem will persist from a particular software package to the more general problem of 
determining starting solutions.}. Here we address this problem by using a minimax approach to
variational elliptic PDE proposed by Zhou and co-workers \cite{LiZhou,LiZhou1} by employing the 
code \cite{Zhou} to generate starting solutions in \texttt{pde2path}, {i.e.}, instead of programming the
continuation method ``by-hand'' \cite{ChenZhou} we glue the available combination of FEM, 
continuation and variational optimization to minimize the programming effort. The goal is to
demonstrate that already today many computational tasks could be simplified substantially if
better interfaces for gluing code would be available, and new code would be designed with this 
problem in mind.\medskip 

To cross-validate the gluing approach, we use the Lane-Emden-Fowler equation, which is a standard example exhibiting 
multiple (unstable) solutions, which are computable by the \texttt{minimax} method \cite{LiZhou,Zhou}. We add a linear 
term term, whose rate is used as the main continuation parameter. Starting solutions are 
first determined by the \texttt{minimax} approach and then continued in \texttt{pde2path}. As a second approach,
the homogeneous zero solution branch is continued numerically, and branch switching is 
performed at bifurcation points to calculate branches, which each contain a single solution
to the Lane-Emden-Fowler equation without the added linear term.\medskip

As the first main new problem for this paper, we consider the Lane-Emden-Fowler equation
with a forcing term. The motivation for modifying the elliptic PDE with a certain forcing term is 
developed by pointing out the analogy to a similar term appearing as a localized microforce 
in the generalized Ginzburg-Landau/Allen-Cahn equation \cite{AllenCahn,AransonKramer,Gurtin,Landau},
where the nonlinearity has the opposite sign of the Lane-Emden-Fowler equation. Although the theory 
of the two classes is quite different due to the different signs of operators \cite{ChenZhouNi}, the 
algebraic structure is similar and the motivation of external microforces also applies to fluid
dynamics applications of the Lane-Emden-Fowler equation. As before, 
we use a \texttt{minimax} approach to determine a set of starting solutions. In this case, there is 
no homogeneous solution branch available for direct numerical continuation, although a two-parameter homotopy 
approach could be tried to locate starting solutions\footnote{This approach becomes from a practical implementation
viewpoint more involved, which is exactly what we aim to avoid here.}. We show that there
is a quite intricate deformation of localized structures along isolated bifurcation curves. 
The bifurcation curves are C-shaped curves, which deform onto thinner structures in the $L^\I$-norm 
as the complexity of the patterns is increased. In particular, the numerical study in this paper illustrates the effect, 
how symmetry-breaking may occur when asymmetric localized forces are added to elliptic PDEs on symmetric 
domains. From an analytical viewpoint, special solutions of nonlinear elliptic PDEs on symmetric domains 
are a very active research area, see {e.g.}~\cite{BartschDAprilePistoia}. The same remark applies to the singularly-perturbed
elliptic PDE related to Lin-Ni-Takagi type problems \cite{AoWeiZeng,LinNiTakagi,LinNiWei}, which has the same
signs for the operators as our modified Lane-Emden-Fowler equation without the asymmetric microforce. Hence, it
is critical to understand, for applications and theory, the effect how inhomogeneous bifurcation branches 
continue globally under microforce loading.\medskip

As a third case study, we consider the 
Caginalp system for phase transitions, such as melting-solidification processes. The domain is chosen as a more
complex, wing-shaped geometry. The motivation for this setup arises directly from (de-)icing problems in
aerospace engineering. As for the previous elliptic PDEs, there are multiple solutions for a fixed parameter 
value, which we compute using \texttt{pde2path} by guessing the right starting solution. 
The \texttt{minimax} approach using the code \cite{Zhou} does not modify in a straightforward way 
to detect the solutions. Therefore, 
a more elaborate programming effort would be required to use the variational starting point, which is 
precisely what we want to avoid here. However, from various recent works, {e.g.}~\cite{XieYuanZhou}, 
one does expect that it
is possible to modify the approach to the Caginalp system, as well as other variational
problems to find starting solutions for continuation; this approach will be considered in 
future work. For the Caginalp system we find an interesting bifurcation diagram with 
potential applications to tune parameter values to access different parts of the so-called
mushy regime between melting and solidification.\medskip

In the last part of the paper, we provide a very brief overview of various software
developing issues for gluing codes in numerical continuation analysis of PDEs. The idea
is to just record, which aspects have been important to facilitate the simplest possible
combination of various algorithms used in this work. The second goal is to provide a 
more detailed outlook, from the viewpoint of gluing computation and an end-user 
perspective of applied nonlinear dynamics, to help developers improve future software
releases. It is hoped that this may help to facilitate the exchange of differential 
equations software between different research communities.\medskip

The paper is structured as follows: In Section \ref{sec:review} we review the basic numerical
tools very briefly and fix the notation. The cross-validation test problem is considered in
Section \ref{sec:test}. The bifurcation analysis for the Lane-Emden-Fowler equation with a
linear term and a localized microforce asymmetric load is carried out in Section 
\ref{sec:microforce}. Numerical accuracy and performance results for the 
algorithms used in this example, are recorded in Appendix \ref{ap:performance}. 
The Caginalp system on an airplane-shaped geometry is numerically 
studied in Section \ref{sec:Caginalp}. Section \ref{sec:software} contains a summary of the
main software issues encountered while carrying out the computations. We conclude with a short
summary in Section \ref{sec:summary}.  

\section{Basic Review of Methods}
\label{sec:review}

In this section, we shall just provide a very quick review of the computational tools we 
are going to use in the gluing-type numerical continuation calculations for elliptic PDEs. Consider
a bounded domain $\Omega\subset \R^2$, $x=(x_1,x_2)^\top\in\Omega$ and $u=u(x)\in\R$. We shall 
assume that the boundary $\partial \Omega$ is piecewise smooth throughout this paper. In this 
paper, we shall restrict\footnote{The idea, and most of the software components, would also 
work for systems of elliptic PDEs. However, it is the goal of this paper to show, how simple 
gluing computation can be, not how difficult it can be made.} to elliptic PDEs 
\cite{Evans,GilbargTrudinger} of the form
\be
\label{eq:ellipticPDE}
-c\Delta u(x) +au(x)-f(u(x),x,\mu)=0,\qquad x\in \Omega\subset \R^2
\ee
where $a\in\R$, $c\in(\R-\{0\})$, $f:\R^3\ra \R$ is sufficiently smooth, $\mu\in\R$ is the main
bifurcation parameter, and 
$\Delta=\frac{\partial^2}{\partial x_1^2}+\frac{\partial^2}{\partial x_2^2}$ is the usual
Laplacian. Furthermore, \eqref{eq:ellipticPDE} will be augmented with either zero 
Dirichlet or Neumann boundary conditions, {i.e.}, we consider
\be
\label{eq:BC}
u(x)=0\quad \text{for $x\in\partial \Omega$}\qquad \text{or}\qquad 
({\vec n}\cdot\nabla u)(x)=0\quad \text{for $x\in\partial \Omega$,} 
\ee
where ${\vec n}\in\R^2$ denotes the outer unit normal vector to $\Omega$, 
$\nabla u=\left(\frac{\partial u}{\partial x_1},\frac{\partial u}{\partial x_2}\right)^\top\in\R^2$
is the gradient of $u$, and `$\cdot$' denotes the standard Euclidean inner product.

\subsection{FEM for Elliptic PDE \& pdetoolbox}
\label{ssec:FEM}

The main underlying PDE solver used for \eqref{eq:ellipticPDE}, via the numerical continuation
software \texttt{pde2path}, is the MatLab \texttt{pdetoolbox} \cite{MatLab2013b_base,MatLabPDE}. A standard 
finite-element method (FEM) is used in the \texttt{pdetoolbox}, which we very briefly recap here to fix some
notation in the case of Dirichlet boundary conditions. Let $H^1_0=H^1_0(\Omega)$ denote the 
usual Sobolev space
\be
H^1_0=\left\{v\in L^2(\Omega):\nabla v\in L^2(\Omega),v|_{\partial \Omega}= 0\right\}.
\ee
Let $\cT=\{T_i\}_{i=1}^{n_t}$ be a triangulation of $\Omega$ consisting of $n_t$ triangles $T_i$
and let $\cN=\{p_k\}_{k=1}^{n_p}$ denote the set of $n_p$ associated node/vertex points $p_k\in \R^2$. 
Of course, $n_p$ can potentially be very large if a fine mesh is needed to correctly resolve
solutions. Then define the space of piecewise-linear finite elements
\be
V_h:=\left\{v\in C(\Omega):\text{$v|_{T_i}$ is affine linear for all $T_i\in\cT$}\right\},
\ee
where the subscript $h$ usually indicates the size of the triangles, {e.g.}, $h$ can be 
taken as the maximum edge length.\medskip 

The weak formulation of \eqref{eq:ellipticPDE}, with a test 
function $v_h\in V_h$, is  
\be
\label{eq:elliptic_weak}
\int_\Omega \left[c(\nabla u\cdot \nabla v_h)(x) +au(x) v_h(x)\right]\txtd x
=\int_\Omega v_h(x) f(u(x),x,\mu)~\txtd x,\quad \forall v_h\in V_h.
\ee
We want to approximate the true solution $u$ by its representation $u_h$ in $V_h$ given by
\be
\label{eq:rep_base}
u_h(x)=\sum_{k=1}^{n_p}U_k\phi_k(x)
\ee
using the nodal basis functions $\phi_k\in V_h$, which are defined via the condition 
$\phi_k(p_j)=\delta_{k,j}$, where $\delta_{k,j}$ is the usual Kronecker delta with 
$\delta_{k,k}=1$ and $\delta_{k,j}=0$ if $k\neq j$. Inserting $u_h\approx u$ into 
\eqref{eq:elliptic_weak} and using $v_h=\phi_k$ as test functions for each $k$, 
leads to a nonlinear system 
\be
\label{eq:ellPDE_discrete}
KU=F(U,\mu) \qquad \text{or}\qquad KU-F(U,\mu)=0,  
\ee
where $U=\{U_k\}_{k=1}^{n_p}\in \R^{n_p}$ is the vector of unknowns and 
$K\in \R^{n_p\times n_p}$, $F:\R^{n_p}\times \R\ra \R^{n_p}$ are assembled automatically in 
the \texttt{pdetoolbox}, {i.e.}, this is the key step, which can be very different for various 
methods as well as different types of PDEs.\medskip 

In addition to mesh generation and assembly of 
the linear system, the \texttt{pdetoolbox} also provides an adaptive mesh refinement algorithm 
\cite{JohnsonEriksson} based on the error estimator
\be
\label{eq:error_estimate}
E_i=E(T_i)=\alpha \|h(f-au_h)\|_{L^2(T_i)}
+\beta c\left(\frac{1}{2} \sum_{\tau\in \partial T_i}h_\tau^2\left[{\vec n_\tau}\cdot \nabla u_h\right]^2\right)^{1/2},
\ee
where ${\vec n_\tau}$ is the outer unit normal and $h_\tau$ the length of the edge $\tau$, 
$h=h(x)$ is the local mesh size for the element $T_i$, the term $\left[{\vec n_\tau}\cdot \nabla u_h\right]$ 
denotes the flux across the element edge, and $\alpha,\beta$ are fixed constants.
In practice, one of the most efficient schemes is to equidistribute the error and split 
those triangles, and relevant adjacent ones, with the highest errors during a mesh refinement step.
The error equidistribution principle has been used for a long time in numerical continuation 
problems for BVPs of ODEs motivated by the COLSYS code \cite{AscherChristiansenRussell}, which 
pre-dates the development of AUTO \cite{Doedel_AUTO2007}. Of course, we are again faced with a
similar challenge as in the discretization step for PDEs as there is a huge variety of ideas and 
estimators for mesh adaptation, not just the solution \eqref{eq:error_estimate} implemented in
the \texttt{pdetoolbox}. 

\subsection{Numerical Continuation \& pde2path}
\label{ssec:continuation}

Here we briefly review the numerical continuation algorithm used in \texttt{pde2path} 
\cite{UeckerWetzelRademacher}; for more theoretical background see \cite{Govaerts,Kuznetsov}. 
For numerical continuation, one views the problem \eqref{eq:ellipticPDE} as a mapping
\be
\label{eq:ellipticPDE_cont}
g(u,\mu):=-c\Delta u(x) +au(x)-f(u(x),x,\mu)=0,\qquad g:\cH\times \R\ra \cH, 
\ee
where we assume, for simplicity, that $\cH$ is a Hilbert space, {e.g.}, $\cH=H^1_0(\Omega)$. The
goal is to trace out a curve (or branch) $z(s):=(u(s),\mu(s))\in \cH\times \R$ of solutions to 
$g(u,\mu)=0$, where $s$ is arclength. However, working in the abstract setting does not lead
immediately to an implementable numerical algorithm and one has to consider the problem via
a spatial discretization as in \eqref{eq:ellPDE_discrete} so one studies
\be
\label{eq:ellipticPDE_cont1}
G(U,\mu):=KU-F(U)=0,\qquad G:\R^{n_p}\times \R\ra \R^{n_p}, 
\ee
where the Hilbert space is just a, potentially very large, Euclidean space $\R^{n_p}$. The goal is
still to trace out a curve $Z(s):=(U(s),\mu(s))\in \R^{n_p}\times \R$. To make $s$ an approximation
to arclength, one augments \eqref{eq:ellipticPDE_cont1} by an extra condition and studies 
\be
H(U,\mu):=\left(G(U,\mu),p(U,\mu,s)\right)^\top=0,
\ee
where $p:\R^{n_p}\times \R\times \R\ra \R$ will be defined below. Suppose we know a point 
$(U(s_0),\mu(s_0))=:(U_0,\mu_0)$ on the curve $Z=Z(s)$ at $s=s_0$ from some initial guess or 
otherwise\footnote{Recall from Section \ref{sec:intro} that this is already a crucial assumption; 
we are going to find initial solutions in Section \ref{ssec:variational} using a specialized 
algorithm.}. Note that we may compute the tangent vector 
\be
\dot{Z}_0:=(\dot{U}_0,\dot{\mu}_0)^\top:=
\left(\left.\frac{\txtd}{\txtd s}U(s)\right|_{s=s_0},\left.\frac{\txtd}{\txtd s}\mu(s)\right|_{s=s_0}
\right)^\top \in\R^{n_p+1}
\ee
to $Z(s)$ at $s=s_0$ by differentiating $G(U(s),\mu(s))=0$ with respect to $s$ at $s_0$ so that
\be
\label{eq:solve_tangent}
(\txtD G)|_{s=s_0}\dot{Z}_0=0,
\ee
where $(\txtD G)|_{s=s_0}$ is the total derivative of $G$ at $s_0$, {i.e.}, finding $\dot{Z}_0$
reduces to solving the linear system \eqref{eq:solve_tangent}. Then, define
\be
p(U,\mu,s):=\xi ~\dot{U}_0\cdot (U(s)-U_0)+(1-\xi)~\dot{\mu}_0(\mu(s)-\mu_0)-(s-s_0),
\ee
where $\xi$ is an algorithm parameter with $0<\xi<1$. In \texttt{pde2path}, $\xi$ is a key parameter 
to ``tune'' or ``control'' the behavior of the numerical continuation algorithm and it is
assumed that $\dot{Z}_0$ is normalized to one in the weighted $\xi$-norm derived from the
inner product
\be
\left\langle (U,\mu)^\top,(W,\nu)^\top\right\rangle_\xi:=\xi~ U\cdot W+(1-\xi)~\mu\nu, 
\qquad \|\dot{Z}_0\|=\sqrt{\langle\dot{Z}_0,\dot{Z}_0\rangle_\xi}\stackrel{!}{=}1.
\ee
Then $p(U,\mu,s)=0$ defines a hyperplane perpendicular, with respect to the inner product 
$\langle\cdot,\cdot\rangle_\xi$, to $\dot{Z}_0$ at a distance $\delta s:=s-s_0$ from 
$(U_0,\mu_0)$. Now, one can finally define the main predictor-corrector steps to compute 
a new point
\be 
(U_1,\mu_1)^\top:=(U(s_1),\mu(s_1))^\top=(U(s_0+\delta s),\mu(s_0+\delta s))^\top.
\ee
The prediction step is to set $(U^1,\mu^1)^\top:=(U_0,\mu_0)^\top+\delta s~\dot{Z}_0$.
One possibility to carry out the correction step is to use a straightforward Newton
iteration
\be
\left(
\begin{array}{c}
U^{j+1}\\
\mu^{j+1}
\end{array}
\right)=
\left(
\begin{array}{c}
U^{j}\\
\mu^{j}
\end{array}
\right)-
\left.\left(
\begin{array}{cc}
\txtD_uG & \txtD_\mu G\\
\xi \dot{U}_0 & (1-\xi)\dot{\mu}_0\\
\end{array}
\right)\right|_{(U^j,\mu^j)}H(U^j,\mu^j)
\ee
Terminating Newton's method at some step $j=J_1$ by a suitable error criterion, we set 
$(U_1,\mu_1):=(U^{J_1},\mu^{J_1})$. So the key algorithm control parameters for the \texttt{pde2path} 
implementation of numerical continuation are the step size $\delta s$, the criterion for terminating 
Newton's method and the parameter $\xi$, which basically balances vertical versus horizontal 
preference along bifurcation curves \cite{UeckerWetzelRademacher}.\medskip 

Of course, there are many 
different variations one may use, {e.g.}, approximating the Jacobian in Newton's method in
different ways \cite{Dijkstraetal}, use Moore-Penrose continuation \cite{Kuznetsov} or replace 
the Newton method by, a sometimes more robust, secant approach. Another issue are bifurcation
points \cite{Kuznetsov}, where one also wants an algorithm to switch between different 
bifurcation branches; we refer the reader for details of these algorithms to \cite{Kuznetsov,Govaerts}
as this is not the main problem considered in this paper.   

\subsection{Multiple Starting Solutions \& minimax}
\label{ssec:variational}

As discussed in Section \ref{sec:intro} and mentioned in Section \ref{ssec:continuation}, one 
main challenge in numerical continuation is to find at least one starting solution. Even if
one solution can be guessed, or an algorithm does converge to a solution from a rough initial
guess, this is often not sufficient. Indeed, frequently not all bifurcation branches connect
to a single branch on which the single known solution lies, or it is very difficult to trace out
all solutions from a single branch using standard branch switching approaches 
\cite{Doedel_AUTO2007}. For example, isolated bifurcation branches, so-called isolas 
\cite{AvitabileDesrochesRodrigues,Keller4}, appear frequently in applications \cite{Becketal,
DesrochesKrauskopfOsinga1,KnoblochLloydSandstedeWagenknecht,McCallaSandstede}. 
In fact, in two-dimensional spatial problems one may potentially find isola structures 
more easily \cite{TaylorDawes}, in comparison to the one-dimensional case\footnote{Quoting from
\cite{TaylorDawes} in the context of a Swift-Hohenberg equation: ´\textit{``Isolas appear for some 
non-periodic boundary conditions in one 
spatial dimension but seem to appear generically in two dimensions.''}}.\medskip 

One approach to systematically look for starting solutions for elliptic PDEs of 
the form \eqref{eq:ellipticPDE} is to consider the associated energy variational functional.
For zero Dirichlet boundary conditions, consider the space $\cH=H_0^1$ and the functional 
\be
\label{eq:functional}
J(u):=\int_\Omega\left[\frac{c}{2} \nabla u(x)\cdot \nabla u(x)+\frac{a}{2}(u(x))^2-F(u(x),x,\mu)\right]\txtd x,
\ee 
where $F$ is the anti-derivative of $f$ with respect to the $u$-component, {i.e.}, 
\be
F(u(x),x,\mu):=\int_0^{u(x)} f(w,x,\mu)~\txtd w.
\ee
When dealing with the energy functional $J(u)$, we shall always assume that the bifurcation 
parameter $\mu$, as well as the other parameters $a,c$, are fixed so we suppress the dependence
in the notation; in particular, we are going to search for a set of starting solutions
for fixed parameter values.\medskip

Let $J'(u)\in \cH$ denote the Fr\'echet derivative of $J$ at $u$. Solutions $u$ of the Euler-Lagrange 
equation $J'(u)=0$ are also called critical points. Recall that 
\be
\int_\Omega J'(u)v~\txtd x=\lim_{\epsilon \ra 0}\frac{J(u+\epsilon v)-J(u)}{\epsilon}=
\int_\Omega g(u,\mu)v~\txtd x
\ee 
where $g(u,\mu)=0$ is the original elliptic PDE, {i.e.}, critical points are weak solutions and 
vice versa \cite{Evans}. By elliptic regularity \cite{Evans}, weak solutions are smooth 
classical solutions. In general, critical points can be maxima, minima
or saddle points. A critical point is non-degenerate if the second Fr\'echet derivative $J''(u)$
exists and is invertible. The Morse index (MI) of a non-degenerate critical point is defined as the maximal 
dimension of a subspace on which $J''(u)$ is negative definite\footnote{Of course, this situation can 
nicely be imagined geometrically for the finite-dimensional case \cite{JostRG}.}. If $\textnormal{MI}=0$, 
then we have a stable minimum for the gradient flow of the energy. For $\textnormal{MI}>0$ there is at 
least one unstable direction. \medskip

Note that minima are relatively easy to find by steepest descent but that saddle points are 
more challenging, particularly if one wants to find multiple ones for fixed parameters. 
Li and Zhou \cite{LiZhou,LiZhou1}, motivated by the work in \cite{ChoiMcKenna,DingCostaChen}, 
developed an algorithm to compute multiple saddle solutions. We shall only outline the basic
proposed strategy here; for details see \cite{LiZhou,LiZhou1}. Let $\tilde{\cH}\subset \cH$
be a subspace of the Hilbert space $\cH$ and consider the unit sphere 
$S_{\tilde{\cH}}:=\{v\in \tilde{\cH}:\|v\|_{\cH}=1\}$. Let $\cL$ be a closed subspace in $\cH$
with orthogonal complement $\cL^\perp$. For each $v\in S_{\cL^\perp}$, define the closed 
half-space
\be
[\cL,v]:=\{tv+w:w\in\cL,t\geq 0\}.
\ee
A set-valued map $P:S_{\cL^\perp}\ra 2^\cH$ is called the peak mapping of $J$ with respect to
$\cL$ if for any $v\in S_{\cL^\perp}$, $P(v)$ is the set of all local maximum points of $J$ in 
$[\cL,v]$. Essentially, the map $P$ collects local maxima in a half-space. A single-valued 
map $p:S_{\cL^\perp}\ra \cH$ is a peak selection of $J$ with respect to $\cL$ if 
\be
p(v)\in P(v),\qquad \forall v\in S_{\cL^\perp}.
\ee   
Basically $p$ selects a single maximum parametrized by $v$; the notion of peak selection can be localized 
by intersecting the relevant sets with neighbourhood of $v$ \cite{LiZhou,LiZhou1}. The main idea
of finding critical points is to restrict to suitable solution submanifolds
\be
\cM:=\{p(v):v\in S_{\cL^\perp}\}.
\ee
Essentially, one would like to think of the manifold $\cM$ as (an approximation to) 
the stable manifold of a saddle critical point.
Then one may use a descent method on $\cM$, in combination with a method to stay on $\cM$ during 
the iteration, to find the saddle point. Of course, one has also to avoid convergence to a saddle
point, which has been found previously. The descent process is a minimization problem and the step
to return to $\cM$ is carried out via a maximization problem, leading to a \texttt{minimax}-type algorithm.
The span of the previously found solutions is the space $\cL$, which one tries to avoid at
future iteration steps, {i.e.}, the search essentially works in the complementary half-space $[\cL,v]$.

Note that the algorithm just outlined, can be viewed as a descent method with a constraint to stay inside 
a certain search space by excluding previously found directions. The entire procedure is easily implemented
via FEM as the FEM discretization just reduces the original problem to a discrete finite set of values
of $u$ at certain points. Then one has to actually solve a large, but finite-dimensional, constrained 
optimization problem. A detailed algorithmic description can be found in \cite{LiZhou,LiZhou1}, the 
implementation of the \texttt{minimax} algorithm we use here is \cite{Zhou}, while a number of 
applications of the method are discussed, {e.g.}, in \cite{ChenZhou,ChenZhou1,WangZhou,XieYuanZhou}.\medskip

It is desirable to try to glue the variational \texttt{minimax} approach to a standard continuation
package, which itself is glued to a standard FEM package. This is precisely, what is carried out 
in practice here for several problems using \texttt{pde2path} as the main focal point. For more details on 
the scientific computing challenges involved, we refer to Section \ref{sec:software}.

\section{A Cross-Validation Test Problem}
\label{sec:test}

As a first step, we are going to compare, for a test problem from the class of elliptic PDEs 
\eqref{eq:ellipticPDE}, two different approaches: 

\begin{itemize}
 \item[(I)] Use the \texttt{minimax} approach to calculate at a fixed parameter value $\mu_0=\mu(s_0)$ several 
 starting solutions $u_0^l=u_0^l(x)$, $l\in\{1,2,\ldots,L\}$ for some $L\geq 2$. Then use numerical 
 continuation for each starting solution, generating multiple solution branches 
 $Z^l(s):=(u^l(s),\mu^l(s))$.
 \item[(II)] Start with a single, preferably simple and easy-to-guess, solution at a fixed parameter 
 value $(u^*_0,\mu^*_0)$ and then continue the branch $Z^*(s)=(u^*(s),\mu^*(s))$. Via branch switching
 at bifurcation points, one can then try to recover the starting set of solutions from (I) as points
 on the branch, {i.e.}, $(u^*(s^l),\mu^*(s^l))=(u^l_0,\mu_0)$ for some values $s^l$.   
\end{itemize}

It is helpful to select a test problem, where we may expect that classical continuation ideas using the 
homotopy approach (II) would suffice, but where the variational problem is nontrivial. As a domain, we choose 
a rectangle $\Omega:=(-l_{x_1},l_{x_1})\times (-l_{x_2},l_{x_2})\subset \R^2$, $l_{x_1},l_{x_2}>0$. Consider 
the elliptic PDE for $u=u(x)$ given by
\be
\label{eq:ac_test}
\left\{
\begin{array}{ll}
0=-\Delta u-\mu u-u^3=g(u,\mu),\quad & x\in\Omega,\\
0=u,\quad & x\in\partial \Omega.
\end{array}\right.
\ee
where $\mu\in\R$ is the main bifurcation parameter. For $\mu=0$, and multiplying the entire equation by $-1$, 
the PDE \eqref{eq:ac_test} is also known as the Lane-Emden-Fowler equation \cite{Chandrasekhar,SerrinZou,Wong}. 
In fact, $\mu=0$ is the standard test case from \cite{Zhou} so we definitely can try to carry out the strategy 
(I) with $\mu_0=0$. Furthermore, $u\equiv 0$ is a solution, so setting $(u^*_0,\mu^*_0)=(0,0)$ is a standard
starting point for strategy (II).\medskip 

\begin{figure}[htbp]
\psfrag{lam}{$\mu$}
\psfrag{x}{$x$}
\psfrag{y}{$y$}
\psfrag{norm}{$\|u\|_\I$}
\psfrag{a}{(a)}
\psfrag{c1}{(c1)}
\psfrag{c2}{(c2)}
\psfrag{d}{(d)}
\psfrag{e}{(e)}
\psfrag{f}{(f)}
\psfrag{b1}{(b1)}
\psfrag{b2}{(b2)}
\psfrag{AC}{$0=-\Delta u-\mu u-u^3$}
	\centering
		\includegraphics[width=1\textwidth]{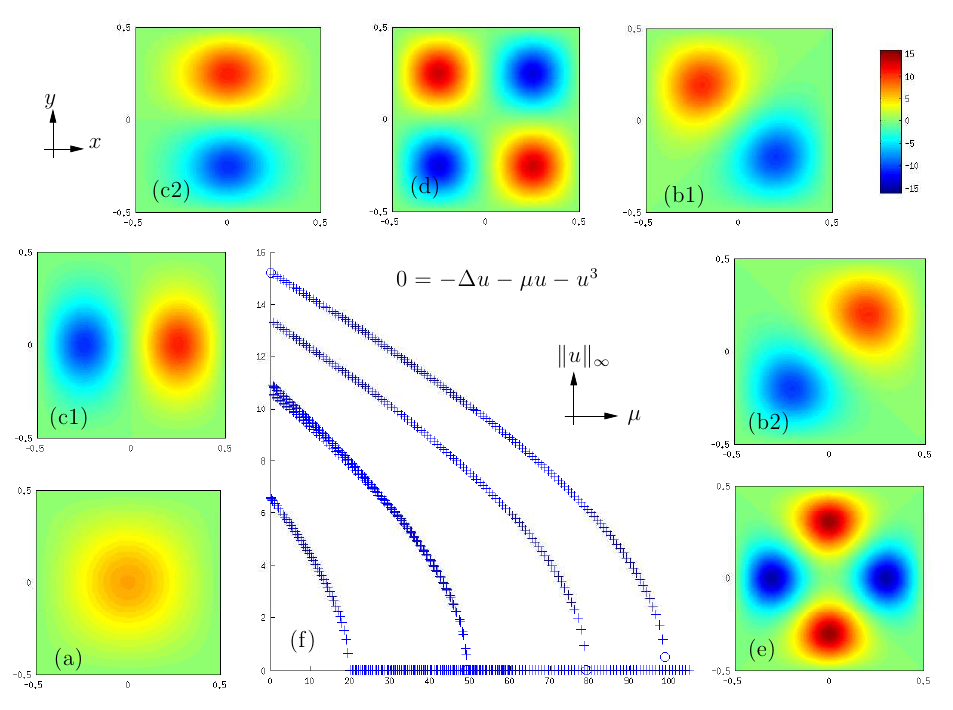}
		\caption{\label{fig:1}Numerical continuation results for \eqref{eq:ac_test} 
		starting from variational solutions found via the \texttt{minimax} algorithm at $\mu=0$. 
		The seven starting solutions are shown in (a)-(e). The four double bump solutions
		(b)-(c) split into two classes with different $L^\I$-norm, where the norm form
		for (b1)-(b2) is smaller than for (c1)-(c2). The norms for (d) and (e) are the two
		largest norms. The part (f) of the figure shows the raw output of \texttt{pde2path} upon 
		continuation; see also the description in the text.}
\end{figure} 

We start with strategy (I). The variational formulation \eqref{eq:functional} of \eqref{eq:ac_test} via 
an energy functional is given, for $u\in H^1_0(\Omega)$ and $\mu_0=0$, by
\be
\label{eq:var_form}
J(u)=\int_\Omega\left[\frac12 \|\nabla u(x)\|^2-\frac14 u^4(x)\right]\txtd x,
\ee 
where $\|\cdot\|$ denotes the usual Euclidean norm in $\R^2$.
The results from the \texttt{minimax} algorithm can be reproduced using \cite{Zhou} and are shown in \ref{fig:1}(a)-(e)
for $l_{x_1}=0.5=l_{x_2}$. Figure \ref{fig:1}(a)-(e) shows four distinct classes of solutions. The 
solution of Morse index $\textnormal{MI}=1$ is shown in Figure \ref{fig:1}(a). Two solutions with two peaks
centered along the coordinate axis of Morse index $\textnormal{MI}=2$ are shown in Figures \ref{fig:1}(c1)-(c2)
and the corresponding two-peak solutions with the peaks along the diagonal and anti-diagonal are displayed in Figures 
\ref{fig:1}(b1)-(b2). Then, there are also two different classes of solutions of Morse index 
$\textnormal{MI}=4$ shown in Figures \ref{fig:1}(d)-(e). Once the starting solutions have been found, they
were saved and re-meshed onto a slightly coarser grid for continuation in \texttt{pde2path}. Figure \ref{fig:1}(f)
shows the resulting bifurcation diagram in $(\mu,\|u\|_\I)$-space, where $\|\cdot\|_\I$ denotes the usual 
$L^\I$-norm. Figure \ref{fig:1}(f) shows the raw output of the \texttt{pde2path} continuation, {i.e.}, the actual 
points computed on the desired seven branches, starting from the seven starting solutions located at $\mu=0$.
The $L^\I$-norms of the seven starting solutions correspond to lexicographical order of the labels, {e.g.}, 
the smallest $L^\I$-norm is in Figure \ref{fig:1}(a), while the largest is in Figure \ref{fig:1}(e).\medskip 

\begin{figure}[htbp]
\psfrag{lam}{$\mu$}
\psfrag{norm}{$\|u\|_\I$}
\psfrag{a}{(a)}
\psfrag{c}{(c)}
\psfrag{d}{(d)}
\psfrag{e}{(e)}
\psfrag{b}{(b)}
	\centering
		\includegraphics[width=0.95\textwidth]{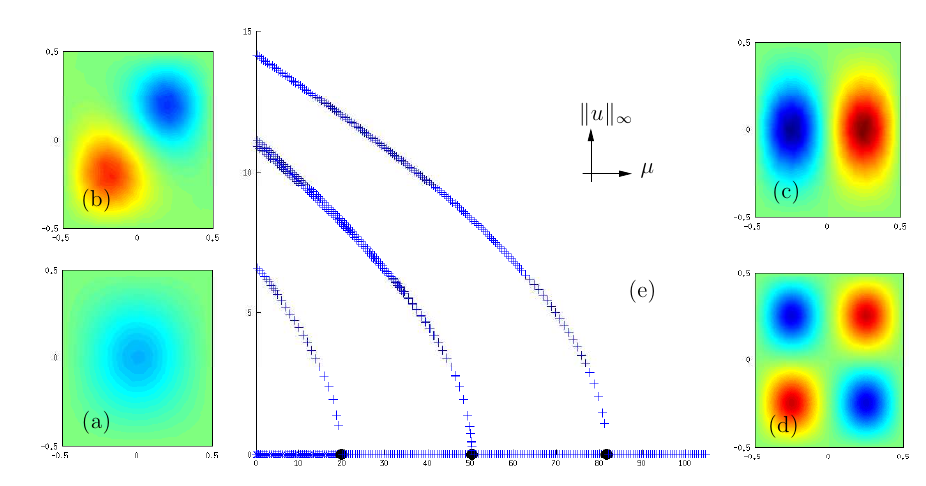}
		\caption{\label{fig:2}Numerical continuation results for \eqref{eq:ac_test} 
		starting from from the zero solution $u\equiv0$ at $\mu=0$; the color coding is
		as in Figure \ref{fig:1}. Further solutions are
		found by branch switching at four detected bifurcation/branch points (marked as black dots).
		The four final ($\mu=0$) solutions on the nontrivial branches are shown in (a)-(d). 
		The two double bump solutions (b)-(c) again have different $L^\I$-norm at $\mu=0$. 
		The branch point for the solution in Figure \ref{fig:1}(e) has not been found in a
		continuation run using the routine \texttt{cont}. However, the branch point can be 
		found when using \texttt{findbif} and is located at $\mu\approx 101.218$. The part 
		(e) of this figure shows the raw output of \texttt{pde2path} upon continuation.}
\end{figure} 

Some of the important algorithmic parameter values for \texttt{pde2path} used during the continuation runs were:
\be
\xi=1/n_p, \quad \texttt{neig}=100, \quad \texttt{dsmax}=0.5, \quad \texttt{hmax}=0.1,
\ee
where \texttt{neig} is the number of computed eigenvalues for stability/bifurcation detection, \texttt{dsmax}
is the maximum step size for $\delta s$ and $\texttt{hmax}$ is the maximum triangle side lengths $h$ for the
(uniform) global triangulation. Otherwise, we used standard parameter values, where it should be noted that
using finite-differencing for the Jacobian seemed to yield a more stable algorithm\footnote{In the forthcoming 
version \texttt{p2p2} of \texttt{pde2path} \cite{DohnalRademacherUeckerWetzel}, this effect, which arises
due to interpolation error, is resolved. I would like to thank Hannes Uecker for explaining me the
details behind this observation.}. The results in Figure
\ref{fig:1}(f) show that all seven initial bifurcation branches connect back to the trivial solution branch
at $(u,\mu)=(0,\mu)$. We also observed that changing some \texttt{pde2path} continuation parameters it
was possible to change the direction the zero branch was continued in once it was reached, either to the
left or to the right; Figure \ref{fig:1} shows all seven branches going to the right, {i.e.}~increasing $\mu$, once
the zero branch was reached. There seems to be a numerical artifact for the branch corresponding to 
the starting solution in Figure \ref{fig:1}(e), where a bifurcation point (blue circle) is detected quite a 
bit before the actual zero branch is reached.\medskip 

For the second approach (II), we continue the trivial zero branch. The results are shown in Figure 
\ref{fig:2}. Four branch points (black dots) were detected when the continuation was run between $\mu=0$ 
and $\mu=110$; again, the last bifurcation point corresponding to the diagonal/anti-diagonal $\textnormal{MI}=4$
solution from Figure \ref{fig:1}(e) was not detected using a standard bifurcation run using \texttt{cont}. However, 
this branch point is correctly detected using the routine \texttt{findbif}\footnote{I would like to thank Hannes 
Uecker for making me aware of this, {i.e.}, that \texttt{findbif} evaluates the eigenvalues, while \texttt{cont} 
checks for a vanishing determinant.}. Branch
switching was performed using \texttt{pde2path} at the bifurcation points. The new non-trivial solutions 
were tracked back up to $\mu=0$, which yields the same - up to symmetry - solutions as recorded for approach (I).
In summary, both approaches essentially lead to the same results and we have definitely cross-validated (I) and (II).   

\section{Adding a Localized Microforce Load}
\label{sec:microforce}

We have observed for the test problem in Section \ref{sec:test} that symmetries lead to difficulties
of detection of certain branches as multiple solution branches come together at a single bifurcation
point. Hence, it is natural to break the symmetry by a certain perturbation. In 
\cite{UeckerWetzelRademacher}, the symmetry is broken for a cubic-quintic Ginzburg-Landau/Allen-Cahn 
type equation of the form 
\be
\label{eq:UWR}
\left\{
\begin{array}{ll}
0=-c\Delta u-\mu u-u^3+u^5,\quad & x\in\Omega,\\
0=u,\quad & x\in\partial \Omega,
\end{array}\right.
\ee
on a rectangle $\Omega:=(-l_{x_1},l_{x_1})\times (-l_{x_2},l_{x_2})\subset \R^2$ with $l_{x_1}=1$ and 
$l_{x_2}=0.9$. Instead of breaking the geometry of the domain, one may also ask, how one expects generic 
symmetry-breaking for Allen-Cahn/Ginzburg-Landau-type PDEs when the equation itself is perturbed. 
One possible mechanism can be 
found in the work of Gurtin \cite{Gurtin}, who derives the Ginzburg-Landau equation from 
basic principles of force balance. The basic motivation for Ginzburg-Landau-type models are
two-phase systems, whose evolution is described on a macroscopic level by an order parameter 
$\rho=\rho(x,t)$ for $(x,t)\in\Omega\times [0,T]$ for some $T>0$. The force balance equations
for a given body lead to the generalized Ginzburg-Landau equation
\be
\beta\frac{\partial \rho}{\partial t}=\nabla \cdot \left[\frac{\partial \Psi}{\partial p}(\rho,\nabla \rho)\right]
-\frac{\partial \Psi}{\partial \rho}(\rho,\nabla \rho)+\gamma, 
\ee
where $p=\nabla \rho$, $\Psi:\R^2\ra \R$ is a given free energy, $\gamma$ is an external
microforce applied to the body and $\beta$ is the constitutive modulus as discussed in \cite{Gurtin}.
Taking a constant positive $\beta$ and the free energy as 
$\Psi(\rho,\nabla \rho)=F(\rho,\mu)+\frac12 \alpha \|\nabla \rho\|^2 $ leads to the Ginzburg-Landau
equation with an applied microforce
\be
\beta\frac{\partial \rho}{\partial t}=\alpha\Delta \rho-\frac{\partial F}{\partial \rho}(\rho,\mu)+\gamma, 
\ee
where $\alpha,\beta$ are positive parameters and $\gamma=\gamma(x)$ is an
external microforce on the body. Frequently, the potential $F$ is chosen as a double-well leading to 
$F'(u)=\mu u-u^3$, which has different signs in comparison to the Lane-Emden-Fowler equation from 
Section \ref{sec:test}. However, the idea that there is a (microscopic) external force acting on 
the problem still seems very reasonable. Hence, we propose to study the elliptic problem  
\be
\label{eq:LEF_force}
\left\{
\begin{array}{ll}
0=-\Delta u-\mu u-u^3+\gamma,\quad & x\in\Omega,\\
0=u,\quad & x\in\partial \Omega,
\end{array}\right.
\ee
where $\gamma=\gamma(x)$ depends in a non-trivial way upon the spatial coordinate and $\Omega:=(-l_{x_1},l_{x_1})\times (-l_{x_2},l_{x_2})\subset \R^2$ with $l_{x_1}=1=l_{x_2}$. 
To break the symmetry encountered in Section \ref{sec:test}, we propose to consider a point-type load
modeled by the function
\be
\label{eq:mforce}
\gamma=\gamma(x):=\left\{
\begin{array}{rl}
\gamma_a\exp\left(-1/\gamma_b[(x_1-\gamma_1)^2+(x_2-\gamma_2)^2)] \right),&\quad\text{if 
$x=(x_1,x_2)^\top \in\Omega$,}\\
0,&\quad \text{if $x\in \partial \Omega$,}
\end{array}
\right.
\ee
and we fix $\gamma_a=10$, $\gamma_b=\frac{1}{10}$, $\gamma_1=0.5=\gamma_2$. Note that 
$\gamma:\Omega\ra [0,+\I)$ is basically localized, with respect to numerical accuracy, in a small 
neighbourhood of the point $(x_1,x_2)=(0.5,0.5)$. Note carefully that using \eqref{eq:mforce} as the forcing 
in \eqref{eq:LEF_force} means that there are no homogeneous constant steady states $u(x)\equiv \text{constant}$ 
available to start the continuation as in Section \ref{sec:test}. This makes \eqref{eq:LEF_force} substantially
more difficult.\medskip

\begin{figure}[htbp]
\psfrag{mu}{$\mu$}
\psfrag{norm}{$\|u\|_\I$}
	\centering
		\includegraphics[width=1\textwidth]{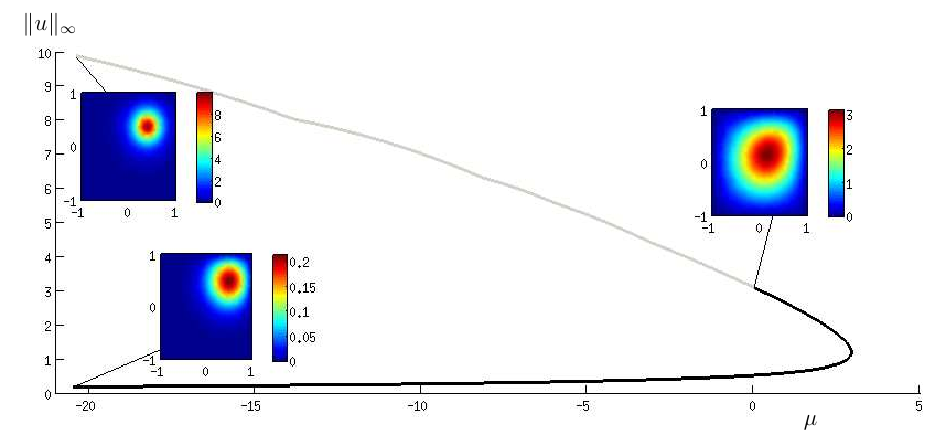}
		\caption{\label{fig:3}Numerical continuation results for \eqref{eq:LEF_force} with 
		external force given by \eqref{eq:mforce}. The starting solution was obtained by 
		the variational \texttt{minimax} approach at $\mu=0$. Then, forward and backward numerical 
		continuation from this starting point leads to two branches (black and grey). The 
		starting solution with Morse index $\textnormal{MI}=1$ is also shown (on the right side), as well
		as the two endpoints on the computed branches (on the left side); for a more detailed description
		see the main text.}
\end{figure} 

\begin{figure}[htbp]
\psfrag{mu}{$\mu$}
\psfrag{norm}{$\|u\|_\I$}
	\centering
		\includegraphics[width=1\textwidth]{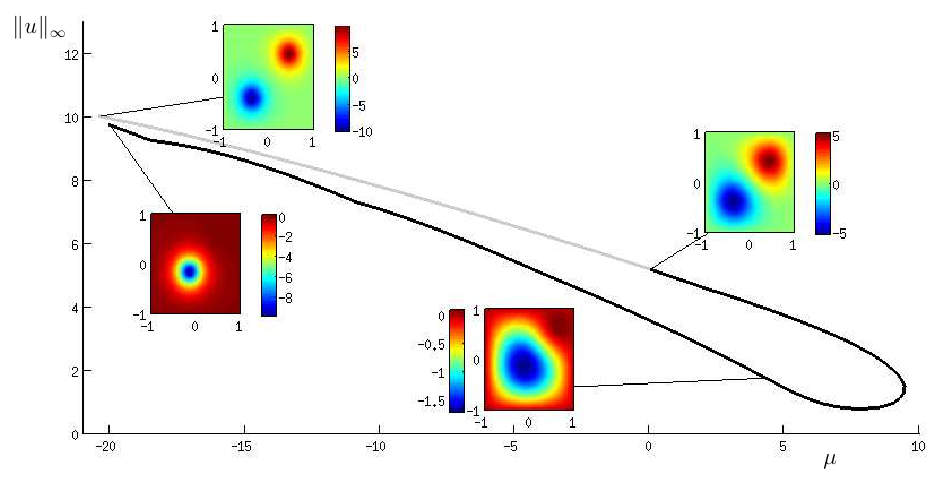}
		\caption{\label{fig:4}Numerical continuation results for \eqref{eq:LEF_force} with 
		external force given by \eqref{eq:mforce}. The starting solution was obtained by 
		the variational \texttt{minimax} approach at $\mu=0$. Then forward and backward numerical 
		continuation from this starting point leads to two branches (black and grey). The 
		starting solution with Morse index $\textnormal{MI}=2$ is also shown; for a more 
		detailed description see the main text.}
\end{figure}

\begin{figure}[htbp]
\psfrag{mu}{$\mu$}
\psfrag{norm}{$\|u\|_\I$}
	\centering
		\includegraphics[width=1\textwidth]{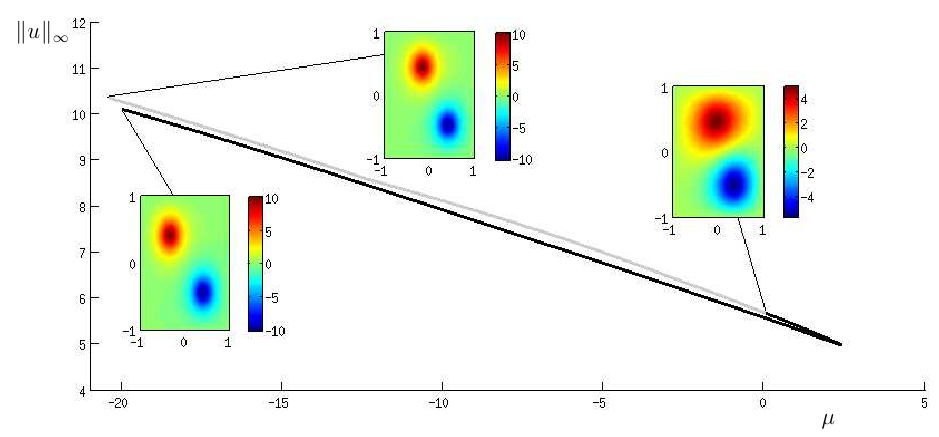}
		\caption{\label{fig:5}Numerical continuation results for \eqref{eq:LEF_force} with 
		external force given by \eqref{eq:mforce}. The starting solution was obtained by 
		the variational \texttt{minimax} approach at $\mu=0$. Then forward and backward numerical 
		continuation from this starting point leads to two branches (black and grey). The 
		starting solution with Morse index $\textnormal{MI}=2$ is also shown. Note that this
		starting solutions is different from the one in Figure \ref{fig:4}; for a more 
		detailed description see the main text.}
\end{figure}

The homotopy continuation strategy (II) with branch switching from Section \ref{sec:test}
is now extremely complicated to carry out. One option is to set $\gamma_a=0$, then start with the zero 
solution, compute the different non-trivial branches using strategy (II), and then continue each branch 
also in the parameter $\gamma_a$ up to $\gamma_a=10$. Albeit certainly possible, it is actually a lot easier 
in this case to follow the gluing strategy (I) as the code can basically be used almost unaltered, just 
by adding $\gamma$ and noticing that the
new energy functional, for $\mu=0$ ({cf.} equation \eqref{eq:var_form}), is just given by 
\be
\label{eq:var_form1}
J(u)=\int_\Omega\left[\frac12 \|\nabla u(x)\|^2-\frac14 u^4(x)+u(x)\gamma(x)\right]\txtd x.
\ee 
The results for the continuation runs are shown in Figures \ref{fig:3}-\ref{fig:5} with 
$\texttt{hmax}=0.07$ and $\texttt{neig}=200$. The results
in Figure \ref{fig:3} show the continuation for the $\textnormal{MI}=1$ solution. The starting solution
at $\mu=0$ has a peak slightly shifted towards the upper right corner of the domain in comparison
to the solution in Figure \ref{fig:1}(a). Continuing forward in $\mu$ leads to the black curve in
Figure \ref{fig:3} along which the solution peak shrinks and moves into the top right corner. On this
curve, a fold point at $\mu=\mu_f\approx 2.92$ occurs, at which the unstable starting solution stabilizes. 
The bottom curve below the fold is expected to be the global attractor for non-steady-state starting solutions 
in the range of for all $\mu\in(-\I,\mu_f]$ for the associated parabolic PDE. 
Running the continuation backwards, {i.e.} with step size $-\delta s$, yields the grey curve. Along
this curve, the solution peak increases and also moves into the top right corner. We also continued both 
parts of the ``(reflected) C-shaped'' branch up to $\mu=-100$ and no branching was detected. In fact, the 
branch shape remains, only the peaks seem to sharpen. In fact, a mesh refinement was helpful for the upper
part of the curve as discussed in Section \ref{ssec:FEM}. The sharpening of the peaks is expected from
a formal scaling argument since we can re-write \eqref{eq:LEF_force} for $\mu\neq 0$ as 
\be
\label{eq:LEF_force1}
\left\{
\begin{array}{ll}
0=-\frac{1}{\mu}\Delta u- u-\frac1\mu u^3+\frac1\mu\gamma,\quad & x\in\Omega,\\
0=u,\quad & x\in\partial \Omega,
\end{array}\right.
\ee
which is a singularly perturbed elliptic PDE as $|\mu| \ra +\I$ (here $\mu\ra -\I$). For 
singularly perturbed elliptic PDEs, concentration phenomena, spike- and/or boundary-layer
solutions are a common phenomenon \cite{LinNiTakagi,Ni,Ni1}. More detailed numerical 
performance calculations were also carried out for this branch and are discussed in Appendix 
\ref{ap:performance}.\medskip 

In Figure \ref{fig:4}, a starting solution from the \texttt{minimax} approach with $\textnormal{MI}=2$
and two peaks along the diagonal is considered. Interestingly, along the lower (black) branch of 
solutions, there is a deformation from a two-peak solution with one minimum and one maximum to a 
single minimum solution. Along the upper (grey) branch, the two-peak structure remains to hold.
Again, we have a C-shaped branch, which folds even tighter back towards itself in the $L^\I$-norm.
Continuing back towards $\mu=-100$ did not lead to the detection of bifurcation points so we 
may actually conjecture that the symmetry-breaking leads to various isolas corresponding to different
peak-structure solutions.\medskip

We also used the \texttt{minimax} method to find another two-peak saddle-solution shown as the
starting solution in Figure \ref{fig:5}. This starting solutions has $\textnormal{MI}=2$
and a slightly broader peak near the top left corner in comparison to the peak in the 
lower right corner. As before, we applied the same forward and backward continuation 
steps and obtain another C-shaped branch, which does not connect to any other solutions
up to $\mu=-100$ and is likely to be isolated from other branches. However, in comparison
to the case of Figure \ref{fig:4}, there is no change in the peak structure along the 
branch and the turning in the $L^\I$-norm is even tighter.\medskip

In summary, the \texttt{minimax} starting strategy turned out to be computationally simple and 
quite effective to determine starting solutions. Furthermore, we may conclude for the 
Lane-Emden-Fowler equation \eqref{eq:LEF_force} with microforce \eqref{eq:mforce} that 
symmetry-breaking by external forces can lead to quite interesting bifurcation branch 
structures, splitting the various branches, which are connected to $u\equiv 0$ when 
$\gamma=0$ as shown in Section \ref{sec:test}. There seem to be many open questions 
regarding the global bifurcation structure for various forces, even for simple 
elliptic problems in the plane. Although
the global structure is easily accessed numerically, one may hope to understand the
problem near a bifurcation point of the trivial branch for $0<|\gamma(x)|\ll1$ analytically
by a perturbation argument and local unfolding; again, this questions seems to be open
for different microforces.  

\section{The Caginalp System}
\label{sec:Caginalp}

In addition to the previous two case studies in Sections \ref{sec:test}-\ref{sec:microforce},
where the \texttt{minimax} starting solutions have been used according to approach (I), one may also
ask for an example where the classical continuation approach is preferable but gluing computation
still plays an important role.\medskip

\begin{figure}[htbp]
\psfrag{x}{$x$}
\psfrag{y}{$y$}
	\centering
		\includegraphics[width=1\textwidth]{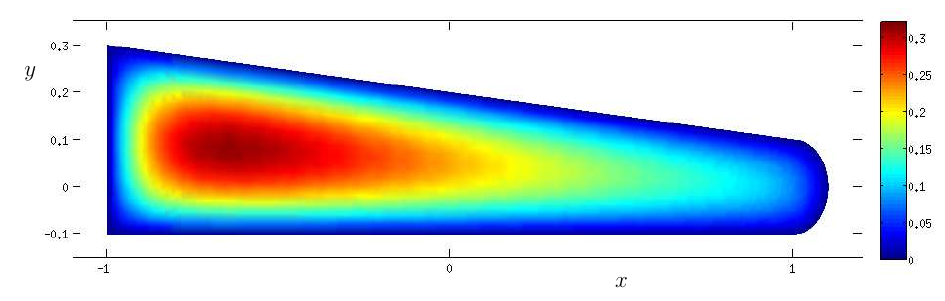}
		\caption{\label{fig:7}Solution of the Poisson equation \eqref{eq:Poisson} computed
		over the basic wing-shaped domain $\Omega$ with a constant point heat source 
		$f\equiv2\cdot 10^3$. This solution is 
		used as an input to the Caginalp system \eqref{eq:Caginalp2}.}
\end{figure}

Consider again a compact domain $\Omega\subset \R^2$, which represents a material or body, and
let $(x,t)\in \Omega\times [0,T]$. The Caginalp system models the interaction between a phase 
field (or order parameter) $\chi=\chi(x,t)$ and the temperature $\vartheta=\vartheta(x,t)$. It
was originally proposed in \cite{Caginalp} and studied using many different techniques 
\cite{GrasselliMiranvilleSchimperna,MiranvilleQuintanilla}. From the modelling perspective, it
describes, {e.g.}, melting-solidification processes between different pure phases. 
Here we use a version of the Caginalp system in the scaling considered in 
\cite{GrasselliPetzeltovaSchimperna} given by
\be
\label{eq:Caginalp}
\left\{
\begin{array}{rcll}
\frac{\partial \vartheta}{\partial t}+ \frac{\partial}{\partial t}\lambda(\chi)  
& = & \Delta \vartheta +f,\quad & \text{if $(x,t)\in\Omega\times [0,T]$,}\\
\frac{\partial \chi}{\partial t} & = & \Delta \chi -(\txtD W)(\chi)+(\txtD\lambda)(\chi)\vartheta,
\quad & \text{if $(x,t)\in\Omega\times [0,T]$,}\\
0&=&\vartheta, \quad & \text{if $(x,t)\in\partial \Omega \times [0,T]$,}\\
0&=&\vec n \cdot \nabla\chi, \quad & \text{if $(x,t)\in\partial \Omega \times [0,T]$,}\\
\vartheta &=&\vartheta_0,\quad & \text{if $(x,t)\in \Omega \times \{0\}$,}\\
\chi &=&\chi_0,\quad & \text{if $(x,t)\in \Omega \times \{0\}$,}
\end{array}
\right.
\ee 
where $W:\R\ra \R$ is a given potential, $f:\Omega \ra \R$ is a heat source/sink and
$\lambda:\R\ra \R$ models the latent heat density. We are just interested in the
stationary Caginalp system under a given fixed heat source. In this case, the system
reduces to
\be
\label{eq:Caginalp1}
\left\{
\begin{array}{rcll}
0  & = & \Delta \vartheta +f,\quad & \text{if $x\in\Omega$,}\\
0 & = & \Delta \chi -(\txtD W)(\chi)+(\txtD\lambda)(\chi)\vartheta,
\quad & \text{if $x\in\Omega$,}\\
0&=&\vartheta, \quad & \text{if $x\in\partial \Omega $,}\\
0&=&\vec n \cdot \nabla\chi, \quad & \text{if $x\in\partial \Omega$,}\\
\end{array}
\right.
\ee 
where all unknown functions $\chi=\chi(x),\vartheta=\vartheta(x)$ only depend upon
the spatial variable. Returning to one modelling purpose of the Caginalp system,
namely melting-solidification processes, we aim to chose a domain $\Omega\subset \R^2$
on which melting and solidification definitely play a crucial role in applications.
One engineering application of interest are airplane wings, where the process of
avoiding ice formation is of critical importance \cite{ThomasCassoniMacArthur}. So
we chose $\Omega$ as a very basic wing-shaped domain as shown in Figure \ref{fig:7}.
\medskip 
 
A first step for the numerical continuation of \eqref{eq:Caginalp1} is to note
that we may use an offline pre-computing step since we have assumed that $f$ is
time-independent, {i.e.}, we may solve a Poisson equation 
\be
\label{eq:Poisson}
\left\{
\begin{array}{rcll}
0  & = & \Delta \vartheta +f,\quad & \text{if $x\in\Omega$,}\\
0&=&\vartheta, \quad & \text{if $x\in\partial \Omega $.}\\
\end{array}
\right.
\ee 
In this context, a gluing framework is very helpful. It is very easy to use the
graphical user-interface provided by the \texttt{pdetoolbox} to define the domain $\Omega$, 
specify $f$, the Dirichlet boundary conditions, the PDE itself, and then solve
the Poisson equation \eqref{eq:Poisson}. The solution for $f\equiv 2000$
is shown in Figure \ref{fig:7}. Hence, if we view $\vartheta(x)=\vartheta$ as a
known temperature input for the Caginalp system \eqref{eq:Caginalp1} then it
remains to study
\be
\label{eq:Caginalp2}
\left\{
\begin{array}{rcll}
0 & = & -\Delta \chi +(\txtD W)(\chi)-(\txtD\lambda)(\chi)\vartheta,
\quad & \text{if $x\in\Omega$,}\\
0&=&\vec n \cdot \nabla\chi, \quad & \text{if $x\in\partial \Omega$,}\\
\end{array}
\right.
\ee    
as an elliptic PDE for the order parameter $\chi$.\medskip

\begin{figure}[htbp]
\psfrag{mu}{$\mu$}
\psfrag{norm}{$\|\chi\|_\I$}
	\centering
		\includegraphics[width=1\textwidth]{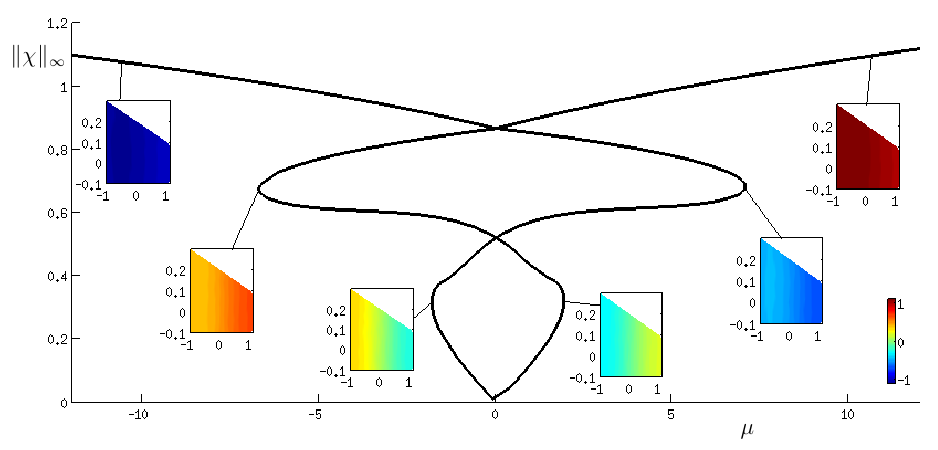}
		\caption{\label{fig:6}Numerical continuation results for \eqref{eq:Caginalp2} with
		$\vartheta=\vartheta(x)$ from Figure \ref{fig:7}, latent heat density \eqref{eq:latentheat}
		and potential \eqref{eq:Ceng}. The starting solution was computed for $\mu=1$ by 
		iterating an initial guess. The remaining part of the diagram was computed using 
		continuation in \texttt{pde2path}. For a more detailed description see the main text.}
\end{figure} 

To study the reduced
Caginalp equation \eqref{eq:Caginalp2} numerically, we still have to specify
the form of the latent heat density $\lambda:\R\ra \R$ and the potential
$W:\R\ra \R$. The latent heat density is frequently modeled and/or empirically determined as 
leading-order affine-linear term with higher-order small-prefactor nonlinear polynomial correction 
terms \cite{RogersYau}. To replicate this modelling, we just consider
\be
\label{eq:latentheat}
\lambda(\tau):=\mu \left(c_0+c_1 \tau +c_2 \tau^2\right),
\ee
where $\mu$ is again the bifurcation parameter and $c_1=1$ and $c_2=0.05$ are
fixed for the computation; note that $c_0$ is not relevant here as only $\txtD \lambda$
appears in \eqref{eq:Caginalp2}. For the potential, we select the standard example  
\be
\label{eq:Ceng}
W(y):=(y^2-1)^2,
\ee
which was already studied by Caginalp \cite{Caginalp}. Note carefully, that
the elliptic PDE \eqref{eq:Caginalp2} has a variational structure with
energy functional
\be
\label{eq:var_form2}
J(\chi)=\int_\Omega\left[\frac12 \|\nabla \chi(x)\|^2+W(\chi)-
\lambda(\chi(x))\vartheta(x)\right]\txtd x.
\ee 
However, the quartic term $\chi^4$ has a positive sign, {i.e.}, the Allen-Cahn/Ginzburg-Landau
sign and not the Lane-Emden-Fowler sign. Still one may formally try to apply a
minimax approach but the \texttt{minimax} code \cite{Zhou} we used for Sections 
\ref{sec:test}-\ref{sec:microforce} does not yield multiple solutions, even upon
minor modifications. Hence, a deeper modification, or a completely different approach
would be required to scan for multiple solutions at a fixed parameter value.\medskip

Naturally, we first have to make sure there really are multiple solutions for some 
fixed parameter value. In this case, the continuation approach can be used to find 
these solutions. In particular, we start at $\mu=1$ by guessing a constant solution.
Under Newton iteration, we actually do reach a true solution of the problem for $\mu=1$.
Then we can use numerical continuation to track the solution branch as shown in Figure
\ref{fig:6}. We observe that there are at least four different solutions at $\mu=1$
so there should definitely be a variational approach to determine them systematically
and we leave this problem as an open problem for future numerical work\footnote{It is
likely that there is a variational algorithm available to search for multiple solutions
in this case but not one, which is freely available for the direct gluing computation
approach proposed here.}.\medskip

The bifurcation diagram in Figure \ref{fig:6} is quite interesting for applications so
we briefly comment on the results. We cut off the bifurcation diagram and only consider
the region with $\|\chi\|_\I\leq 1$ as $\chi=\pm1$ represent the pure states and we are interested 
in the transition (or so-called mushy) regime between the pure states. Suppose we start in the 
bifurcation diagram with an (almost) pure 
state $\chi(x)\approx 1$
for all $x\in\Omega$ and then change $\mu$. This leads to a transition sequence involving 
four folds as shown in Figure \ref{fig:6}. If we interpret $\chi(x)\approx 1$ as a melted
(or liquid) state, then the solution curve at the first fold assumes a mixed/mushy state
where solidification starts non-uniformly over the domain, with a focus near the airplane body.
This focus reverses as the solution branch is traversed further and at the third fold,
when the solidification focuses on the tip. In the last step, the system transitions to the almost fully solid state
after the fourth fold has been traversed. Although we shall not pursue these ideas here further,
it is clear that these observations should have interesting applications in engineering 
applications to build airplane wings and/or during a de-icing process. We leave a detailed 
parametric numerical continuation study for melting/solidification processes for future work.     

\section{Perspectives on Software Development}
\label{sec:software}

Although we have successfully demonstrated that gluing computing between FEM, continuation
and minimax can be efficient and lead to very interesting practical results for elliptic PDEs, 
we have not commented on several
important practical scientific computing issues. In particular, the question is why the computations
here worked in a relatively straightforward way, while gluing computation can become very
complicated in many cases when one tries to track patterns via numerical continuation; the thesis 
\cite{Avitabile} is an excellent example, how challenging such a computational approach can become. 
The following issues seem to be very helpful to keep in mind for further
software development\footnote{Of course, the various issues are written from the viewpoint
of a user in applied nonlinear dynamics, who is interested in gluing computations.}:
\medskip   

\textbf{High-level language:} A key ingredient to efficiently glue different parts
required for the computation is to use a high-level programming language. Here \texttt{MatLab} 
\cite{MatLab2013b_base} is used. However, an excellent non-commercial alternative would be 
\texttt{Python} \cite{Python}, which contains aspects designed for gluing computation and is already
efficiently used in various scientific computing packages for differential equations 
\cite{Doedel_AUTO2007,Hoffmanetal,GuyerWheelerWarren,Cimrman,Herouxetal}. The problem with working directly 
with FEM, continuation or variational-PDE packages via a fast-computation, but lower-level, programming 
language, {e.g.}~\texttt{C}, \texttt{C++} or \texttt{Fortran}, is that even apparently simple-looking 
tasks of transferring data, adding on functionality, problem formulation, interlinking of algorithms, 
cross-validation, testing and visualization, become incredibly complicated and 
time-consuming. To really make a gluing approach work on the basis of a lower-level language, 
one has to fully integrate multiple software packages into a single, necessarily constrained,
environment. Although this may be very desirable in standardized industrial applications, it does
not adequately represent software development and flexibility requirements in an academic environment.
Only a higher-level language provides the required flexibility. Of course, ones has to give up a
slight bit of computational efficiency but overall, this trade-off seems worthwhile.\medskip 

In this context, it should be noted that such an integration of input-output via a high-level language
should probably be made a design principle from the start. If wrapper-functionality is added
later on, it is frequently still necessary to fully comprehend the underlying low-level code
to actually use the wrapper.\medskip 

\textbf{Algorithmic blocks:} Based upon the previous point, one should make sure to design 
self-consistent blocks of code, which interface/communicate directly with a high-level language
but run on a lower-level fast language for computational purposes. For example, in the context 
considered in this paper one might want to split up the process into the following components:

\begin{enumerate}
 \item[(A1)] problem formulation and definitions;
 \item[(A2)] mesh generation and error-estimator based mesh refinement;
 \item[(A3)] discretization of the PDE, {i.e.}, conversion to a nonlinear algebraic system;
 \item[(A4)] algorithms for generating starting solutions sets;
 \item[(A5)] numerical continuation and bifurcation detection; 
 \item[(A6)] efficient, fast numerical linear algebra;
 \item[(A7)] data analysis and visualization.
\end{enumerate}

Of course, not all the different steps are necessary for a given class of PDEs or chosen discretization
algorithm. For example, sometimes one may omit (A4) due to the availability of starting solutions as for approach (II) 
in Section \ref{sec:test}, or use spectral or other mesh-free methods to generate the nonlinear system
of algebraic equations to be solved in (A3).\medskip  

\textbf{Problem description:} Another natural question is, how we should formulate the PDE meshing, discretization
and continuation problem in the high-level gluing language? In this regard, \texttt{pde2path} builds upon the 
successful ideas implemented in \texttt{AUTO}. The main idea is to have only very few files that specify all the problem 
details completely: one (or two) files to specify the elliptic PDE, one file to initialize all the data (domain, initial guess,
algorithmic parameters, {etc.}) and one structure to track the current state of the numerical problem. Those ideas turned out, 
again, to be extremely robust and helpful to carry out the gluing computations in this paper. We also completely
glued the \texttt{minimax}-algorithm \cite{Zhou} to \texttt{pde2path} by using only data from standard \texttt{pde2path} problem 
description files in \texttt{minimax}. In principle, this could be implemented permanently into \texttt{pde2path} but
there is a new version released soon \cite{DohnalRademacherUeckerWetzel}. The variational starting solution approach
will be implemented in this new version in future work\footnote{The new version \texttt{p2p2} is not backward compatible
with \texttt{pde2path} so we postpone this step to future work. It is expected that compatibility is guaranteed from the version \texttt{p2p2} onwards \cite{DohnalRademacherUeckerWetzel}.}.\medskip
 
\textbf{PDE discretization:} By now, there are a large number of different software packages available to automatically
generate meshes, nonlinear equations and adaptive mesh refinements for various classes of PDEs; see Section \ref{sec:intro}.
Many of these software packages would provide excellent tools for understanding dynamics and patterns in a wide
variety of applications a lot better. In fact, the barrier does not lie in the packages themselves but the difficulty
to access the problem formulation. For example, it is relatively rare that one can simply 
provide a problem description and obtain, \emph{in a simple format}, the discretized nonlinear equations as an output; of course, it is
always \emph{possible} but the goal is to make gluing computation \emph{easy}. However, there seems significant progress in this 
direction, {e.g.}~see \cite{Hoffmanetal}. The demonstration we have given here for three problems also strongly points towards the conjecture
that practical gluing computations for dynamical systems analysis of PDEs will become a lot easier in the very near future.\medskip

\textbf{Continuation algorithms:} The main barrier to make continuation algorithms more applicable to spatial
problems was to move beyond the class of two-point BVPs and provide a more generic package suitable for PDEs.
This is one current disadvantage with the \texttt{pde2path} focus we considered in this paper since \texttt{pde2path} links 
internally a continuation algorithm directly to the MatLab \texttt{pdetoolbox}. Therefore, it is not easy in practice 
to replace the steps (A3) and (A5) as they are tied together in \texttt{pde2path}. However, this could be 
resolved by using a generic continuation toolbox, such as the recently developed continuation core \texttt{COCO} 
\cite{DankowiczSchilder} or certain tools for large-dimensional linear systems such as \texttt{LOCA} \cite{Salingeretal} or other recently developed codes \cite{Bindeletal}. The main point is, that in the future one must 
take advantage of the advanced numerical analysis discretization schemes to reduce the computational linear algebra
effort more efficiently. Overall, there is also significant recent progress in this direction supporting the conjecture stated in the last paragraph.

\section{Summary}
\label{sec:summary}

In this paper we have addressed the issue of gluing computation at the interface between numerical
PDEs, continuation methods, scientific computing and variational theory in the context of elliptic PDEs for 
three different equations. We 
have shown that efficient integration of various algorithms within the framework of a high-level
programming framework can be possible. A combination of \texttt{pde2path}, the MatLab \texttt{pdetoolbox} 
and a \texttt{minimax}-algorithm has been glued and cross-validated using the Lane-Emden-Fowler equation 
with a linear term as a test problem. As a second step, we argued
that a natural symmetry-breaking mechanism, which is not based upon the domain geometry, could be localized asymmetric
microscopic forces which also arise in generalized Ginzburg-Landau/Allen-Cahn equations. We found
interesting inverse-C-shaped bifurcation curves. There is numerical evidence to conjecture that these curves are isolas, which could be analytically captured in the regime of a very small microforce. Along the bifurcation curves, complicated deformation of patterns can already take place under a simple near-localized external microforce. Then we 
proposed to study as a third example the Caginalp system for melting-solidification processes
on a wing-shaped geometry. In this case, the variational \texttt{minimax}-algorithm could not directly be glued despite
the existence of multiple solutions at fixed parameters. This illustrated the limitations of the process. For the Caginalp system, we also computed an interesting phase transition diagram for melting-solidification processes, which pointed towards potentially useful aerospace and engineering applications. Finally, we also 
commented on scientific computing issues encountered during the process.\medskip  

\textbf{Acknowledgments:} I would like to thank the Austrian Academy of Sciences ({\"{O}AW})
for support via an APART fellowship. I also acknowledge support of the European Commission 
(EC/REA) via a Marie-Curie International Re-integration Grant. I would like to thank John
Guckenheimer, Hannes Uecker and Mathieu Desroches for very helpful comments on earlier versions of this 
manuscript. Furthermore, I would like to thank two anonymous referees, whose comments have helped to improve 
the paper.


\newpage

\appendix

\section{Numerical Performance Tests}
\label{ap:performance}

In this appendix, we collect some performance data for the numerical continuation calculations.
There are two main reasons, why it seems useful to collect those results here: (R1) There is
to be relatively little detailed data available for \texttt{pde2path} used in
combination with the MatLab \texttt{pdetoolbox} and (R2) we want to verify that the computational
speed is sufficient to allow for direct numerical experiments on a current standard desktop computer
without long waiting times to access the structure of different bifurcation branches.\medskip

\begin{figure}[htbp]
\psfrag{h}{$h$}
\psfrag{a}{(a)}
\psfrag{b}{(b)}
\psfrag{c}{(c)}
\psfrag{np}{$nt$}
\psfrag{nt}{$np$}
\psfrag{sec}{$\tau$}
	\centering
		\includegraphics[width=1\textwidth]{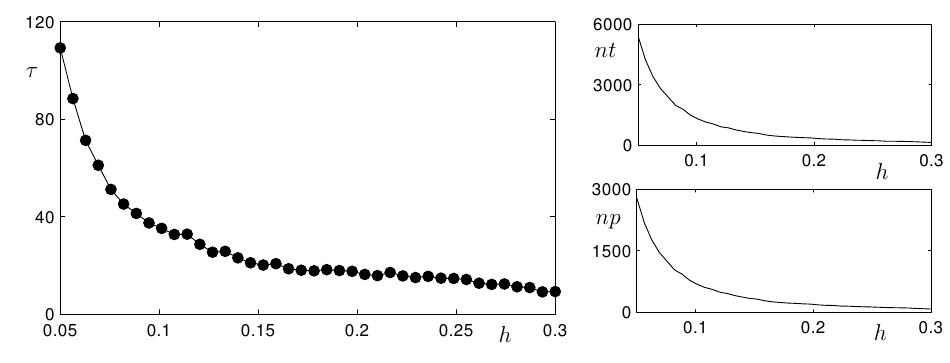}
		\caption{\label{fig:8}Performance results for continuation calculation of the entire curve in
		Figure \ref{fig:3} using different meshes with maximal triangle length $h$. The panel 
		on the left shows the computing time $\tau$ in seconds, while the two panels on the right
		just record the number of points $np$ and number of triangles $nt$ in the mesh depending upon $h$.
		The dots mark the computed points in the left figure and the lines are just linear interpolation.
		In the right part of the figure, the points have been omitted and just the linear interpolation
		is shown.}
\end{figure} 

As a test problem, we re-compute the branch shown in Figure \ref{fig:3} arising from the 
Lane-Emden-Fowler equation with a linear term \eqref{eq:LEF_force} and the external microforce 
\eqref{eq:mforce}. The main reason, why the curve in Figure \ref{fig:3} seems to be a good 
standard test problem, is that it has a fold point, where the solution changes stability, a
single localized solution structure with Morse index one, a singularly-perturbed spike-type
structure on the upper branch and a regular smooth solution structure on the lower branch. 
Hence, despite
the algebraic simplicity of the problem, it is expected that many other solution structures for 
elliptic problems are going to have similar features on bifurcation branches, just with 
higher Morse indices.\medskip

\begin{figure}[htbp]
\psfrag{mu}{$\mu$}
\psfrag{LI}{$\|u\|_\I$}
\psfrag{L2}{$\|u\|_2$}
\psfrag{u10}{\scriptsize u10}
\psfrag{u20}{\scriptsize u20}
\psfrag{u30}{\scriptsize u30}
\psfrag{l10}{\scriptsize l10}
\psfrag{l20}{\scriptsize l20}
\psfrag{l30}{\scriptsize l30}
	\centering
		\includegraphics[width=1\textwidth]{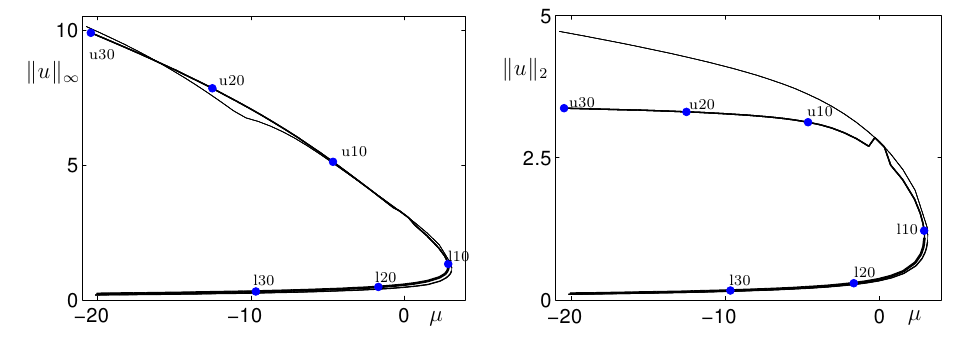}
		\caption{\label{fig:9}Continuation of the entire curve from Figure \ref{fig:3} using a fixed 
		mesh with $h=0.3$ in comparison to an adaptive mesh with parameter \eqref{eq:amesh}. The 
		more accurate branch with adaptive meshing has also six different computed points marked 
		on it (in blue), with associated adaptive meshes as shown in Figure \ref{fig:10}. Note that
		in the two figures (left/right) on the norm of the solution is different on the vertical axis
		{i.e.}, once we use the $L^\I(\Omega)$-norm and once the $L^2(\Omega)$-norm.}
\end{figure} 

\begin{figure}[htbp]
\psfrag{u10}{\small u10}
\psfrag{u20}{\small u20}
\psfrag{u30}{\small u30}
\psfrag{l10}{\small l10}
\psfrag{l20}{\small l20}
\psfrag{l30}{\small l30}
	\centering
		\includegraphics[width=1\textwidth]{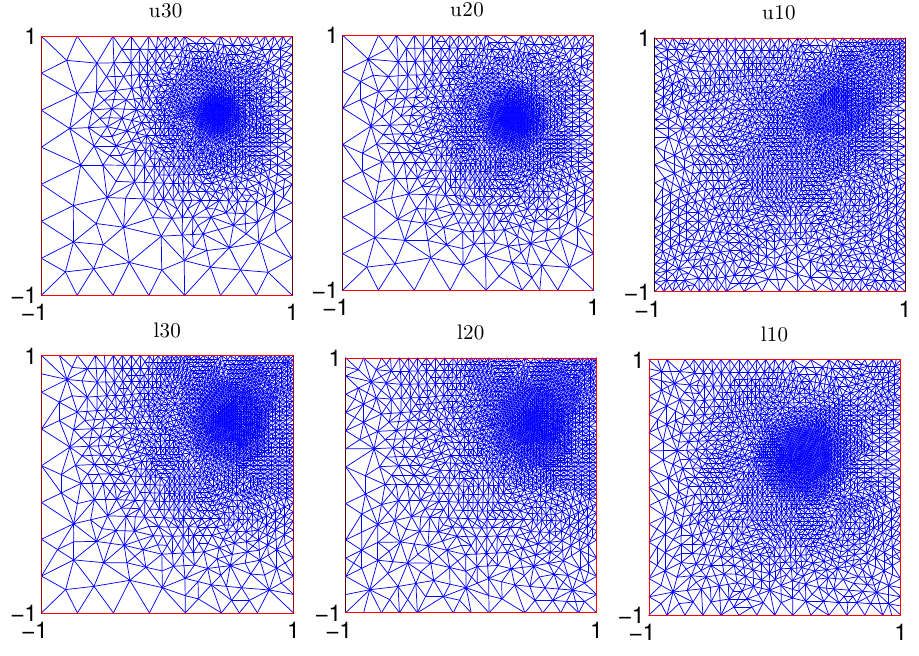}
		\caption{\label{fig:10}Different adaptive meshes corresponding to the points marked 
		in blue in Figure \ref{fig:9}. Observe that 
		the effect of trying to automatically focus 
		on the spike solution for the singularly-perturbed case along the upper branch is clearly 
		visible in the solutions $u20$ and $u30$. On the lower branch, the solution is still localized
		but does not develop a sharp large spike.}
\end{figure}

As a first test, we are interested in the influence of mesh size on speed and accuracy of
computing the entire branch in Figure \ref{fig:3} using \texttt{pde2path}. We use the standard
\texttt{pde2path} algorithm parameters \cite{UeckerWetzelRademacher} and only make the 
following changes 
\be
\texttt{p.ds=0.1},\qquad \texttt{p.tol=1e-6},\qquad \texttt{p.xi=1/p.np},  
\ee
where \texttt{p} is the main dictonary-type continuation structure of \texttt{pde2path} and \texttt{p.np} is the
number of points in the mesh, \texttt{p.xi} is the continuation tuning parameter $\xi$ discussed in
Section \ref{ssec:continuation}, \texttt{p.tol} is the Newton solver tolerance and \texttt{p.ds}
the initial continuation step size. We use a maximum triangle edge length of $h$ and initialize 
the mesh via \texttt{stanmesh} of \texttt{pde2path} respecitively \texttt{initmesh} of the 
\texttt{pdetoolbox}, which yields a Delaunay triangulation of the domain. We run the entire 
continuation, with spatial mesh adaptation turned off (\texttt{p.amod=0}) of the curve shown 
in Figure \ref{fig:3} for different values of $h$. The computation time results are shown in Figure 
\ref{fig:8}. For a quite wide range of mesh sizes the calculation is sufficiently fast to obtain 
the entire bifurcation curve on a standard current desktop computer\footnote{Basic details 
of the desktop computer setup used: Intel Core i5-4430 CPU @ 3.00GHz processor (quadcore), Kingston 16GB system
memory, WDC WD10EZRX-00L 1TB harddrive, Ubuntu 12.04.4 LTS operating system, MatLab R2013a computational 
platform.}.\medskip

\begin{figure}[htbp]
\psfrag{xi}{$\xi$}
\psfrag{neig}{$neig$}
\psfrag{t}{$\tau$}
\psfrag{a}{(a)}
\psfrag{a1}{(a1)}
\psfrag{b}{(b)}
	\centering
		\includegraphics[width=1\textwidth]{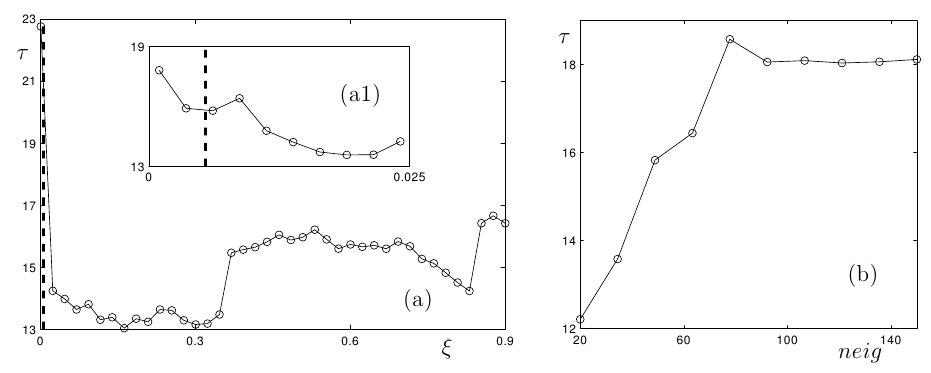}
		\caption{\label{fig:11}Dependence of the continuation calculations on the parameter $\xi$
		(see Section \ref{ssec:continuation}) and the number of eigenvalues calculated. Again the entire
		curve from Figure \ref{fig:3} is calculated for different sets of algorithmic parameters. The circles 
		mark the computed points and the lines are just linear interpolation again. (a) Shows
		a variation of $\xi$ plotted against the computing time $\tau$ (in seconds). The vertical dashed 
		thick line is the suggested value of $\xi$ from \cite{UeckerWetzelRademacher}. The inset (a1)
		shows a finer computation near this suggested value. (b) Computation time depending on the number
		of eigenvalues $neig$ computed to determine stability.}
\end{figure}

Then, it is natural to ask about the accuracy of the calculation and which dynamical features
are captured along the continuation branch in comparison to the mesh size. In Figure \ref{fig:9}
we show the bifurcation curve computed once with a fixed maximum mesh triangle size mesh $h=0.3$
in comparison to the case with mesh adaptation as discussed in Section \ref{ssec:FEM}. The mesh
adaptation was chosen by setting
\be
\label{eq:amesh}
\texttt{p.amod=5},\qquad \texttt{p.ngen=10},\qquad \texttt{p.maxt=4000},  
\ee
in \texttt{pde2path}, which means that mesh adaptation is performed every 5 steps along the
bifurcation branch and one aims for a maximum number of 4000 triangles in at most 10 refinement
steps as discussed in Section \ref{ssec:FEM}. Note that the solution branch in the 
$L^\I(\Omega)$-norm is very well-captured, even on the coarse mesh without mesh adaptation.
However, on the upper part of the bifurcation branch, decreasing $\mu$ leads to a 
singularly-perturbed problem as discussed in Section \ref{sec:microforce} and the 
$L^2(\Omega)$-norm of the two computations clearly differs. The kink on the branch in the 
$L^2(\Omega)$-norm for the lower curve arises from the fact that the continuation is started 
with a coarse mesh ($h=0.3$) and then mesh adaptation only starts to act at the fifth point 
and corrects the solution norm to its actual value.\medskip 

Figure \ref{fig:10} demonstrates
that the mesh is automatically adapted exactly around the forming spike solution as $\mu$
decreases. Therefore, we could also have found the singularly-perturbed nature of the problem
without the scaling  argument leading to equation \eqref{eq:LEF_force1} in an efficient and
immediate way. Note that the branches in Figure \ref{fig:9} are expected to diverge in the 
$L^2(\Omega)$-norm as the coarse mesh solution (upper branch) cannot resolve this norm 
near a spike correctly. However, all indications point to the fact that this is the same 
branch as the maximum, the region where the solution is nonzero as well as the shape of 
the spike are the same. It has also been checked that the Newton solver tolerance does not
seem to affect the continuation run as the same branch was obtained for all 
$\texttt{p.tol}\in(10^{-11},10^{-3})$.\medskip
 
In Figure \ref{fig:11}(a), we investigate the dependence of the continuation time upon the
algorithmic tuning parameter $\xi$ as discussed in Section \ref{ssec:continuation}. The 
vertical dashed line marks the suggested setting taking $\xi$ as the inverse number of points
in the mesh \texttt{p.xi=1/p.np}. Interestingly, this setting does not seem to be the best 
from the computational efficiency point of view as it becomes more difficult to track the branch
around the fold point. However, chosing $\xi$ large is also not good as this makes the
more horizontal parts of the branch slower to track.\medskip 

In Figure \ref{fig:11}(b), the role of the computation time depending upon the number of 
eigenvalues calculated \texttt{p.neig} is considered. Obviously, we expect that the computational
effort increases, when more eigenvalues are requested. However, there seems to be at least
some saturation effect. Hence, we may conclude that there is certainly
potential to optimize continuation runs via algorithmic parameters in the case that larger 
systems, or many solution branches, are to be considered.

\end{document}